\documentclass[a4 paper]
{article}

\usepackage[T2A]{fontenc}
\usepackage[cp1251]{inputenc}
\usepackage[centertags]{amsmath}
\usepackage{amsfonts}
\usepackage{amssymb}
\usepackage{amsthm}
\usepackage{amsmath}
\usepackage{amscd}
\usepackage{latexsym}
\numberwithin{equation}{section}
 \newtheorem
   {theo}{Theorem}
    \newtheorem
   {lem}{Lemma}
\newtheorem {cit}{Theorem}
\title{
Accuracy of approximation  of subharmonic functions by logarithms
of moduli of analytic ones in Chebyshev metrics  }
\author{ Markiyan Girnyk
 }
\date{}
\begin{document}
\maketitle \begin{abstract}It is known that a subharmonic function
of finite order  $\rho$ can be approximated by the logarithm of
the modulus of an entire function at the point $z$ outside an
exceptional set up to  $C\log|z|$. In this article we prove that
if such an approximation is made more precise, i. e. a constant
 $C$ decreases, then, beginning with  $C=\rho/4$, the size of the
exceptional set enlarges substantially. Similar results are proved
for subharmonic functions of infinite order and functions
subharmonic in the unit disk. These theorems improve and
complement a result by Yulmukhametov.
\end{abstract}\large \section{Introduction}
 The results on approximation in various metrics of a function subharmonic
 in a domain  $\Omega$ by the logarithm of the modulus of an
 analytic function have numerous applications
 (see, for instance, [1-6]) since it is often easier to construct
  a subharmonic function rather than  an analytic one with desired asymptotic properties.
Mainly in this connection, a problem arises of approximating
elements of a broader set of subharmonic functions by elements of
a smaller set of the logarithms of moduli of analytic ones. The
problem of that approximation was investigated by a number of
authors. The case, when the Riesz measure is concentrated on
curves, is considered in  [7-9]. Azarin [10] was the first to find
such an approximation in the general form in the class of
functions subharmonic in the plane and having finite order of
growth. This approximation for arbitrary subharmonic functions in
arbitrary domains is considered in the article by Yulmukhametov
[11], where, in particular, the  precise approximation in the
plane $\mathbb{C}$ is found. The research by Yulmukhametov is
extended and developed in the different directions by a number of
authors (see a survey in [6]
and recent works [12-15]).

 We use the principal results and the standard notations of potential theory
 [16]. Let us recall some of them. We denote $D(a, r):=\{z:|z-a|<r\},\,\,
 C(a, r):=\{z:|z-a|\le r\},\,\,
 S(a, r):=\{z:|z-a|=r\},\,\,A(t, T):=\{z:t\le |z|\le T\},\,m_2$-- the plane Lebesgue measure,
 the letters $C$ with indices
 stand for positive constants, in parentheses we indicate the dependence on parameters.
   We denote by $FM$ the set of the subsets $[1,\infty)$ having finite measure and the subsets  $L$ of the interval
    $[0,\,1)$,
    for which $\int_L(1-x)^{-2}\,dx < \infty$. Let $u(z)$ be a subharmonic function in the plane  $\mathbb{C}$ or in the disk  $\mathbb{D}:=D(0,1)$
    , then $\mu_u$ is its Riesz measure,  $B(r,u):=\max\{u(z):z\in C(0,r)\}\,$ is the maximum,
    $n(a,r,u):= \mu(C(a,r)),\,n(r):=n(0,r,u)$ are the counting functions of the Riesz measure. The set of points,
    $(\beta,s)$ -normal with respect to a measure $\mu$ on a domain $\Omega,
    \mathrm{supp}\,\,
    \mu \subset \overline {\Omega}$ , is set by the condition

    $N(\beta,s,\mu):=\{z\in\Omega:\forall t\in (0,s)(\mu (D(z,t))\le\beta t)\},\,\, N(\beta,s,\mu_1,...,\mu_n)
    =\cap_{j=1}^{j=n}N(\beta,s,\mu_j).$

     The notation $a\asymp b$ means $|a|\le \mathrm{Const}|b|$ and
      $|b|\le
    \mathrm{Const}|a|$.
 The content of our work is deeply connected with two results
 by Yulmukhametov, which we cite
 (in somewhat modified but equivalent formulation)[11, pp. 278-282].
\begin{cit}Let  $u(z)\,$ be a subharmonic function of finite order
$\rho$ and a number $\alpha > \rho$ .

 Then there exist an entire function $f(z)$, a constant $C_1=a_0+a_1\alpha,\,a_1 > 1,$
depending only on $\alpha$, and an exceptional set $E$, depending
on the functions  $u(z),\,f(z)$, and the number $\alpha$, such
that
\begin{equation}
    \left|u(z)-\log |f(z)|\right|\le C_1(\alpha) \log|z|
    ,\, z\notin E,\end{equation}
 with $E\subset\cup_j D(z_j,r_j) $ and
\begin{equation} \sum_{R\le|z_j|<2R }r_j=o(R^{\rho -\alpha}),\quad R \to \infty.
     \end{equation}
\end{cit}

\begin{cit}Let an entire function $f(z)$ satisfy the relation
\begin{equation}||z|-\log|f(z)||=o(\log |z|),\, E \not \ni z \to
\infty.
 \end{equation}
Then for every covering $E$ by disks $\{D(z_j,r_j)\}$ with
uniformly bounded radii  and every number $\varepsilon>0$ the
inequality
\begin{equation} \sum_{R\le|z_j|<2R }r_j \ge
R^{1-\varepsilon},\,R>R(\varepsilon).
\end{equation}holds.
\end{cit}

The aim of our work is to improve and to make more precise theorem
B and, also, to prove similar results  for subharmonic functions
of infinite order and functions subharmonic in the unit disk.
 Let us formulate our results.
     \begin{theo}Let a number $\rho>0$ and an entire function $f$ satisfy the inequality
          \begin{equation}||z|^\rho-\log|f(z)||\le C_2\log
          |z|,\,z \notin E,
     \end{equation}
      where a set $E$ is contained in the union of disks  $\{D(z_j,r_j)\}$ with $r_j \le |z_j|^{1-\rho /2 +\varepsilon},$
      \newline $\,\varepsilon >0$.
      Then the estimate
           \begin{equation} \sum_{R\le|z_j|<2R }r_j \ge
R^{1+\rho /2-2C_2-3\varepsilon},\,R>R(\varepsilon),
     \end{equation} holds.
   \end{theo}
 \begin{theo} Let a subharmonic function
   $u(z)$ be of infinite order and a number $\varepsilon >0$.
 Then there exist an entire function  $f(z)$, a constant
   $ C_2(\varepsilon)$, and a set $E$
  such that for all $z\notin E$ the inequality
\begin{equation}
    \left|u(z)-\log |f(z)|\right|\le C_3(\varepsilon)\left( \log|z|+\log B(|z|,u)    \right)\end{equation}
 holds.
The exceptional set $E$, depending on the functions $u(z),\, f(z)$
and a number
   $\varepsilon$,is contained in the union of disks
   $D(z_J, r_j),\,j\in \mathbb{N}$,
with
\begin{equation}
\sum\limits_{R\le |z_j|<R+(\log B(R,u))^{-1}
}r_j=o({B}(R,u)^{-\varepsilon}),\, R\to \infty,\,R\notin L \in FM.
    \end{equation}
 \end{theo}

 \begin{theo}
Let $v(R)\in C^2[0,\infty)$ be a convex function of $\log R$,
$v(R)\to \infty $, when $R\to \infty$, and the conditions
 $$
 \min \{\log v^{(k)}(x):x\in [R-2/\log v(R),R+2/\log v(R)]\}\sim
  $$
\begin{equation}\sim
\max \{\log v^{(k)}(x):x\in [R-2/\log v(R),R+2/\log v(R)]\}\sim
\log v^{(k)}(R),\,\,R\to \infty,\, k=0,1,2,
  \end{equation}hold.
Let an entire function $f(z)$ satisfy the inequality
\begin{equation}
|v(|z|)-\log|f(z)||< C_4 \log v(|z|),\,z\notin E,
  \end{equation}
 where a set $E$ is contained in the union of disks
  $D(z_j, r_j),\,j\in \mathbb{N}$ ,
  such that $r_j<v(|z_j|)^{-1-\varepsilon},\,\varepsilon >0$.
Then
  \begin{equation}
\sum\limits_{R\le |z_j|<R+(\log B(R,u))^{-1} }r_j \ge
v(R)^{1/2-2C_4-2\varepsilon},\,R>R(\varepsilon).
  \end{equation}
 \end{theo}

 \begin{theo}Let a number $\rho>0$ and an analytic function in the disk  $f(z)$
 satisfy the inequality
          \begin{equation}||1-z|^{-\rho}-\log|f(z)||\le C_5\log
          \frac 1
         { |1-z|},\,z \notin E,
     \end{equation}
      where a set $E$ is contained in the union of disks $\{D(z_j,r_j)\}$, with $r_j \le |1-z_j|^{1+\rho /2 },\,$.
      Then the estimate
           \begin{equation} \sum_{R\le|z_j|<R+(1-R)/2 }r_j \ge
(1-R)^{-\rho /2+2C_5+3\varepsilon},\,R(\varepsilon)<R<1,
     \end{equation}holds.
   \end{theo}
 \begin{theo} Let a function $u(z)$ be subharmonic in $\mathbb{D}$
  and of infinite order and a number $\varepsilon >0$.
  Then there exist a function $f(z)$ analytic in $\mathbb{D}$ , a constant
   $ C_6(\varepsilon)$, and a set  $E$ such that for all
  $z\notin E$ the inequality
\begin{equation}
    \left|u(z)-\log |f(z)|\right|\le C_6(\varepsilon)\log B(|z|,u),\end{equation}
holds.
 The exceptional set $E$, depending on the functions $u(z),\, f(z)$ and a number
   $\varepsilon$, is contained in the union of disks
   $D(z_J, r_j),\,j\in \mathbb{N}$,
  with
\begin{equation}
\sum\limits_{R\le |z_j|<R+(1-R)^2(\log B(R,u))^{-1}
}r_j=o({B}(R,u)^{-\varepsilon}),\, R\to 1,\,R\notin L \in FM.
    \end{equation}
 \end{theo}

 \begin{theo}
Let $v(R)\in C^2[0,1)$ be a convex function of  $\log R$, $v(R)\to
\infty $, when $R\to 1$,  and the conditions
 $$
 \min \{\log v^{(k)}(x):x\in [R-2(1-R)^2/\log v(R),R+2(1-R)^2/\log v(R)]\}\sim
  $$
$$ \sim \max \{\log v^{(k)}(x):x\in [R-2(1-R)^2/\log
v(R),R+2(1-R)^2/\log v(R)]\}\sim $$\begin{equation}\sim  \log
v^{(k)}(R),\,\,R\to 1,\, k=0,1,2,
  \end{equation} hold.
Let a function $f(z)$ analytic in $\mathbb{D}$ satisfy the
inequality
\begin{equation}
|v(|z|)-\log|f(z)||< C_7 \log v(|z|),\,z\notin E,
  \end{equation}
 and a set $E$ is contained in the union of disks
   $D(z_j, r_j),\,j\in \mathbb{N}$ ,
  such that $r_j<v(|z_j|)^{-1/2+\varepsilon},\,\varepsilon >0$.
  Then

  \begin{equation}
\sum\limits_{R\le |z_j|<R+(1-R)^2(\log B(R,u))^{-1} }r_j \ge
v(R)^{1/2-2C_7-2\varepsilon},\,\,R(\varepsilon)<R<1.
  \end{equation}
 \end{theo}

 Let us comment these statements. Theorem 1 sharpens Theorem В,
announced by Yulmukhametov. The restriction from above on the
radii of the covering of exceptional set in that and other similar
theorems is necessary: the disks
\newline $\{D((-1)^n 2^{n-1},2^n)\}_{n\in \mathbb{N}}$ cover all the plane
, but $$\sum_{|z_j|<R }r_j\le 2R. $$ Upper bound for the radii of
the disks of the exceptional set covering in Theorems 1, 3, 4, 6 is
greater than or equal to the radii of the disks in Theorem A and
similar theorems respectively, and, therefore, it seems natural. If
one considers the covering of the plane by the disks of the form
$D(z,|z|^{1-\rho/2})$, with multiplicity at most 3, then for that
covering  $$\sum_{R\le|z_j|<2R }r_j \asymp R^{1+\rho/2} ,$$ i.e. a
set $E$ can be exceptional. It is easy to see that functions,
satisfying the conditions of Theorems
 3 and 6, exist, for instance,
$v(R):=\exp_k R,v(R):=\exp_k (1/(1- R)) $, where $\exp_k$ is the
$k$-th iteration of exponent. Theorem 2 sharpens the result in
[17]. Theorems 4-6 complement theorems of [18].

The cases of the plane, the disk, and the functions of finite and
infinite orders are of specific character, because of that we have
not succeeded to represent Theorems
 1, 3, 4, 6 as the implications of one theorem, although the proofs of them are based on the same idea
.

 I am thankful to R. Yulmukhametov, who informed me the idea of the unpublished proof of Theorem B.
 I express my gratitude to V. Eiderman, who recommended to use
  the maximum in place of the Blumental growth function,
and that enables us to improve the result in
 [17]. The participants of Lviv seminar on complex analysis, particularly I. Chyzhykov and O. Skaskiv
 , carefully discussed the proofs and made a few useful comments,
 and my pleasant duty is to notice their work thankfully.

 \section{ Proofs of
  results}

 For the reader's convenience we cite a few theorems used in what
 follows.
\begin{cit} $\mathop
    {\rm[11, c. 275]}$ Let $u(z)$ be a function subharmonic in an
    unbounded domain
  $\Omega$. Then there exists a function $f(z)$,
 analytic in $\Omega$, which for all the $ (\beta,s)$ --normal points with respect to the Riesz
 measures
$\mu_u$ and $\mu_{\log |f|}$ satisfies the inequality
 \begin{equation}
 |u(z)-\log|f(z)||\le \left|
 \int_{D(0,2)}\log|z-\xi|\,d\mu_u(\xi)\right|+C_8 |\log s|+C_9
 \log|z|+\beta s(|\log s|+1)+C_{10}.
 \end{equation}

 Constants $C_8,\,C_9,\,C_{10}$ do not depend on the domain $\Omega$ and the functions
  $u(z),\,f(z)$.

\end{cit}
\begin{cit}

$\mathop
    {\rm[11, c. 275]}$ Let $u(z)$ be a function subharmonic in a
    bounded domain
  $\Omega$ . Then there exists a function $f(z)$,
  analytic in $\Omega$ , which for all the $ (\beta,s)$-normal points with respect to the Riesz
 measures  $\mu_u$ and $\mu_{\log |f|}$ satisfies the inequality
 \begin{equation}
 |u(z)-\log|f(z)||\le
 C_{11} |\log s|+C_{12}|\log \mathop
    {\rm diam}\Omega|
+\beta s(|\log s|+1)+C_{13}.
 \end{equation}

 Constants $C_{11},\,C_{12},\,C_{13 }$ do not depend on the domain $\Omega$ and the functions
  $u(z),\,f(z)$.

\end{cit}

 As Yulmukhametov notes, the upper bounds for the function  $\log|f(z)|$ , obtained in Theorems C and D
 , hold for  the $ (\beta,s)$ -normal points with respect only to the measure
  $\mu_u$.
\begin{cit}$\mathop
    {\rm[19, p. 121]}$
Let  $v(r)\to \infty$, when $r\to \infty$, be a continuous
nondecreasing functions defined on $[r_0,\infty)$ a number $\delta
>0$. Then for all  $r\ge r_0 $ , except, possibly, a set of finite measure
 $L \in FM$ , the inequality
\begin{equation}
v\left(r+\frac 2 {\log v(r)}\right)<v(r) ^{1+\delta}
 \end{equation}holds.

\end{cit}
The following theorem is easy obtained from Theorem E by the
change of variable.
\begin{cit}Let  $v(r)\to \infty$, when $r\to 1$, be a continuous
nondecreasing functions defined on $[r_0,1)$ a number $\delta
>0$
 . Then for all  $r_0 \le
r < 1 $ , except, possibly, a set
 $L \in FM$ , the inequality
\begin{equation}
v\left(r+\frac {2(1-r)^2} {\log v(r)}\right)<v(r) ^{1+\delta}
 \end{equation} holds.
\end{cit}

At first, we shall prove an elementary but important lemma.
\begin{lem}Let $v(x)$ be a twice continuously differentable
function, defined on $[0,\infty)$ or on  $[0,1)$, and
 $$
 h(x,v):=v(T)\frac{\log x -\log R}{\log T-\log R}+v(R)\frac{\log T -\log x}{\log T-\log
 R},
 $$
 where numbers $R,\,T,\,R<T,$ belong to the range of definition of
 the
 function
 $v$ . Then for every
$x \in [R,\,T]$ the inequality
 $$|h(x,v) -v(x)|\le C_{14}(T-R)^2(v'((R+T)/2)/R+\max_{\xi \in
 [R,T]}|v''(\xi)|)
 $$ holds.

\end{lem}

\textbf{The proof of Lemma 1.} By the Taylor formula in the
Lagrange form we have
\begin{equation}
v(x)=v\left(\frac{R+T}{2}\right)+v'\left(\frac{R+T}{2}\right)\left(x-\frac{R+T}{2}\right)+\frac{v''(\xi)}
{2 }\left(x-\frac{R+T}{2}\right )^2,
 \end{equation}
 where $\xi \in [R,T]$,
$$ \frac{\log (1+\frac{x-R}{R})}{\log
(1+\frac{T-R}{R})}=\frac{\frac {x-R}{R} + O((\frac {x-R}{ R})^2)
}{\frac {T-R}{R} + O((\frac {T-R }{R})^2) }=\frac{x-R + O(\frac
{(x-R)^2}{ R}) }{ {T-R} + O(\frac {(T-R)^2 }{R}) }=
\frac{x-R}{T-R}(1+O((T-R)/R)),$$
\begin{equation}(T-R)/R \to 0,
\end{equation}
\begin{equation}
\frac{\log(1+\frac{T-R}{R})-\log(1+\frac{x-R}{R})}{\log(1+\frac{T-R}{R})}=\frac{T-x}{T-R}(1+O((T-R)/R)),\,\,(T-R)/R
\to 0.
\end{equation}

Taking into account (2.5)-(2.7), we can write down $$h(x,v)
-v(x)=(v(T)-v(x))\frac{\log (1+\frac{x-R}{R})}{\log
(1+\frac{T-R}{R})}+(v(R)-v(x))\frac{\log(1+\frac{T-R}{R})-\log(1+\frac{x-R}{R})}{\log(1+\frac{T-R}{R})}=
 $$
 $$=
 \left(v\left(\frac{T+R}{2}\right)+v'\left(\frac{T+R}{2}\right)\left(T-\frac{T+R}{2}\right)+
 \frac {v''(\xi_1)}{2}\left(T-\frac{T+R}{2}\right)^2-\right.$$ $$\left.-v\left(\frac{T+R}{2}\right)-
 v'\left(\frac{T+R}{2}\right)\left(x-\frac{T+R}{2}\right)-\frac
 {v''(\xi_2)}{2}\left(x-\frac{T+R}{2}\right)^2\right)\times$$
 $$\times\frac{x-R}{T-R}(1+O((T-R)/R))+
 $$
 $$
 +\left(v\left(\frac{T+R}{2}\right)+v'\left(\frac{T+R}{2}\right)\left(R-\frac{T+R}{2}\right)+
 \frac {v''(\xi_3)}{2}\left(R-\frac{T+R}{2}\right)^2-\right.$$
 $$\left.
 -v\left(\frac{T+R}{2}\right)-
 v'\left(\frac{T+R}{2}\right)\left(x-\frac{T+R}{2}\right)
 -\frac {v''\left(\xi_4\right)}{2}\left(x-\frac{T+R}{2}\right)^2\right
 )\times$$
 $$\times
 \frac
 {T-x}{T-R}(1+O((T-R)/R))=
 $$
 $$
 =
 v'\left(\frac{T+R}{2}\right)O((T-R)^2/R)+\left(\frac {v''(\xi_1)}{2}-\frac {v''(\xi_2)}{2}+\frac
 {v''(\xi_3)}{2}-\frac {v''(\xi_4)}{2}\right)O((T-R)^2), $$
and from this the statement of Lemma 1 follows.

 \textbf{The proof of Theorem 1.}
 We start from the exposition of the proof idea. We shall first prove that every disk
 of the form  $D(a,\,|a|^{1-\rho/2})$, where $f(a)=0$, contains a quite large exceptional set. Next, it
will be proved that every disk with somewhat greater radius has
the same property without demand on the center of the disk to be
zero of the function $f$. To finish the proof, we shall put
sufficiently many nonoverlapping  disks with enlarged radius into
the annulus
 $A(R,2R)$. Comparing the areas of the annulus and the disks, we shall obtain estimate
 (1.6).
\begin{lem}Let $h(z,v,A)$ be the minimal harmonic majorant of the
subharmonic function $v(z):=|z|^\rho$, where $\rho
> 0$, in the annulus $A:=A(R-R^{1-\rho/2},\,R+R^{1-\rho/2})$. Then
 $$
 \max \{ h(z,v,A)-v(z):z\in A(R-R^{1-\rho/2},\,R+R^{1-\rho/2})\}\le C_{15},
  $$
 where a constant $C_{15}$ does not depend on $R$.
\end{lem}
\textbf{The proof of Lemma 2.}
 It is sufficient to consider the
difference  $h(x,v,A)-v(x)$ , as $$
 h(z,v,A):=v(R+R^{1-\rho/2})\frac{\log
|z|- \log (R-R^{1-\rho/2}) }{\log (R+R^{1-\rho/2}) -\log
R-R^{1-\rho/2})}+ $$ $$ + v(R-R^{1-\rho/2})\frac{\log |z|-\log
(R+R^{1-\rho/2}) }{\log ( R-R^{1-\rho/2})-\log (R+R^{1-\rho/2})}
$$
 is a radial symmetric function. We apply Lemma 1 with
$v(x)=x^\rho$. Since $\max \{ |v'(x)|: x \in
[R-R^{1-\rho/2},R+R^{1-\rho/2}]\}\asymp R^{\rho-1}$, and $\max \{
|v''(x)|: x \in [R-R^{1-\rho/2},R+R^{1-\rho/2}]\}\asymp
R^{\rho-2}$, we obtain the statement of Lemma 2 .

\begin{lem} Let the conditions  of Theorem 1 hold and $f(a)=0$. Then for every
covering of the intersection of $E$ and the disk
$D(a,\,|a|^{1-\rho/2})$  by the disks the sum of their  radii is
not smaller than
 $|a|^{ 1- \rho/2- 2 C_2 -\varepsilon}.$
\end{lem}

\textbf{The proof of Lemma 3.}
 Let us consider the circumferences
$S(a,\,t|a|^{1-\rho/2}),\,t\in [1/2,\,3/4]$ . If all of them
intersect the set  $E$, then the intersection of $E$ and the
disk $D:=D(a,\,|a|^{1-\rho/2})$ can be covered by disks only under
the condition that the sum of their radii is not less than
 $ \frac{1}{8}|a|^{1-\rho/2}$ : compare the radial projections of
 $E$ and the exceptional disks on radius  $D$ .
 If there exists a circumference
$S:=S(a,p|a|^{1-\rho/2}),\,p\in [1/2,\,3/4]$, which does not
intersect $E$, then inside this circumference the estimate
\begin{equation}|h(z,\log|f|,D)-\log|f(z)||<(2C_{2}+\varepsilon)\log|z|,\,\,z
\notin E,
\end{equation} is true,
 where  $h(z,\log |f|,D)$ is the minimal harmonic majorant of
the subharmonic function  $\log |f|$ in the disk  $D$.

Let us prove (2.8). We denote by $h(z,v,D)$ the minimal harmonic
majorant of the subharmonic function  $v(z)=|z|^\rho$ in the disk
$D$ . By the definition of the minimal harmonic majorant and Lemma
2
\begin{equation} 0
 \le h(z,v,D)-v(z)\le h(z,v,A)-v(z)\le C_{15}.\end{equation}
From the maximum principle for harmonic functions it follows that
in $D$
\begin{equation}|h(z,v,D)-h(z,\log |f|,D)|\le |h(z_0,v,D)-h(z_0,\log |f|,D)|
\end{equation} is valid,
 where $z_0 \in S$. Taking into account the conjecture on $S$ , we extend  estimate
 (2.10)to
\begin{equation} |h(z_0,v,D)-h(z_0,\log
|f|,D)|=|v(z_0)-\log|f(z_0)||< C_2 \log |z_0|<
(C_2+\varepsilon)\log |z|.
\end{equation}
     Statement (2.8) follows from (1.5), (2.9)-(2.11).

On the other hand, by the Poisson-Jensen formula for the function
$\log |f|$ in the disk $D$
 $$\log|f(z)|=h(z,\log|f|,D)-\sum\limits_{a_n \in
  D}g(z,a_n,D),
 $$ where $g(z,a_n,D)$ is the Green function of the domain $D$ having the pole at the zero $a_n$ of the function $f$
 ,we have
$$ |h(z,\log|f|,D)-\log|f(z)||=\left|\sum\limits_{a_n \in
D}g(z,a_n,D)\right| =\sum\limits_{a_n \in
 D}g(z,a_n,D).$$
 Applying known properties of Green functions, from this we obtain
\begin{equation}
|h(z,\log|f|, D)-\log|f(z)|| \ge g(z,a,D)= \log \frac
{p|a|^{1-\rho/2}}{|z-a|}.
\end{equation}
 From (2.12) it follows that the disk $D(a,|a|^{1-\rho/2-2C_2-\varepsilon}) \subset
  E$.
 Having proved Lemmas 2 and 3, let us extend the proof of Theorem 1.
 Let us consider an arbitrary number $b \in \mathbb{C}$ and the disk
  $D(b,|b|^{1-\rho/2+\varepsilon})$ with $\varepsilon >0$. There are
  two possible cases: i) the disk $D(b,|b|^{1-\rho/2-2C_2-\varepsilon})\subset
  E$, and that is all right; ii) the difference $D(b,|b|^{1-\rho/2-2C_2-\varepsilon})\backslash
  E \neq \emptyset$, and in this case we replace $b$ by any number
  $c$
  in the difference.
 We again have the following alternative
  : all the circumferences $S(c,t|c|^{1-\rho/2+\varepsilon}),\,t \in [1/2,3/4],$ intersect  $E$, or there
  exists a circumference
  $S(c,p|c|^{1-\rho/2+\varepsilon}),\,p \in [1/2,3/4]$, which does not intersect  $E$. In the first case, for every disk covering
  the sum of the radii is not less than $\frac {1}{8}|c|^{1-\rho/2+\varepsilon}\asymp \frac {1}{8}|b|^{1-\rho/2+\varepsilon}$ when $|b|$
  is enough large.
  For the second case, in the disk
  $D(c,\frac{3}{4}|b|^{1-\rho/2+\varepsilon})$ there exists a number $a$ such that $f(a)=0$. Let us justify this. Let $E\cap S(c, p|c|^{1-\rho/2+\varepsilon})= \varnothing$ and $f$
  has no zeros in the disk
   $D(c,p|c|^{1-\rho/2+\varepsilon})$, hence, $\log|f|$
 is a harmonic function there. We can write down  the Poisson-Jensen  formula
$$
\frac{1}{2\pi}\int\limits_0^{2\pi}(\log|f(c+p|c|^{1-\rho/2+\varepsilon}e^{i\varphi})|-
|c+p|c|^{1-\rho/2+\varepsilon}e^{i\varphi}|^\rho)\,d\varphi=$$
\begin{equation}\log|f(c)|-|c|^{1-\rho/2+\varepsilon}+\int\limits_0^{p|c|^{1-\rho/2+\varepsilon}}
\frac{n(c,t,\log|f|)-n(c,t,v)}{t}\,dt.\end{equation}

 From the definition of the set $E$ it follows that the right-hand side in  (2.13) is $O(\log |c|)=O(\log |b|)$. On the other hand,
 the Riesz measure  $n(c,t,\log |f|)=0$ for every $t \in [0,p|c|^{1-\rho/2+\varepsilon}]$.
 The Riesz measure of the function
   $v(z)=|z|^{\rho}$
  satisfies the relation
   $$
   d\mu_v(z)=\frac{1}{2\pi}\Delta v
   \,dm_2(z)=\frac{1}{2\pi}\rho^2|z|^{\rho-2}\,dm_2(z),
   $$
    where $m_2(z)$ is the plane Lebesgue measure. Hence, we obtain that
   $n(c,t,v)\asymp |c|^{\rho-2}t^2\asymp |b|^{\rho-2}t^2,\, 0\le t\le p|c|^{1-\rho/2+\varepsilon}.$
  Thus, the right-hand side in (2.13) has order $\asymp |b|^{\rho-2}|b|^{2-\rho+2\varepsilon}=|b|^{2\varepsilon}$,
  and it is a contradiction. Therefore, the disk $D(c,\frac {3}{4}|c|^{1-\rho/2+\varepsilon})$ contains
   a zero $a$ of the function $f$.
  The disk $D(a,|a|^{1-\rho/2})\subset D(b,|b|^{1-\rho/2+\varepsilon})$
  for sufficiently large  $|b|$ . Having applied Lemma 3, we conclude that for any disk covering of the intersection
   $E\cap D(b,|b|^{1-\rho/2+\varepsilon})$ the sum of their radii is not less than
  $|b|^{1-\rho/2-2C_2-\varepsilon}$.

  From the comparison of the areas it follows that in the annulus $ A=A(R,2R)$
   the quantity
   $\asymp R^2 /  R^{2-\rho +2\varepsilon}=R^{\rho-2\varepsilon}$ of  the nonoverlapping disks of the
   form $D(b,|b|^{1-\rho/2+\varepsilon})$ can be placed, with
   every disk  containing such a portion of the exceptional
   set that for every disk covering of it the sum of  radii of
   those disks is not less than
   $\asymp R^{1-\rho/2-2C_2-\varepsilon}$.
   Because the radii of the disk covering $E$ are bounded by a number $\asymp
   R^{1-\rho/2+\varepsilon}$, and by this reason  are counted at
   most finite times, with multiplicity not depending on
   $R$, then the total sum of their radii is greater than or equal to
    $\asymp R^{1+\rho/2-2C_2-3\varepsilon}.$

   \textbf{The proof of Theorem
    2.}
   By Theorem С, under the conditions of Theorem 2 there exists an entire function
    $f$, satisfying  (1.7) for all the points $z \in N(\beta, s, \mu_u, \mu_{\log|f|})$
   , where $\beta= B(r,u)^{1+\varepsilon}, s= \beta^{-1}$.

  We turn to the proof of estimate (1.8) of the size of the exceptional set
    $E= \mathbb{C}\setminus  N(\beta, s, \mu_u, \mu_{\log|f|}) $.
    We put $\mu:=\mu_u+\mu_{\log|f|}$. Later the inequality
   $(0< \delta < \varepsilon)$
\begin{equation}n(r,\mu)<C_{16}B(r,u)^{1+\delta},\, r \notin L \in
FM,
\end{equation} will be proved.
  Under the assumption that  (2.14) is true, we estimate the size of $E$. By the definition of a normal point
  with respect to a measure, every point
     $z \in E$
   is the center of a disk $D(z,r_z)$ such that $\mu (D(z,r_z)) \ge B(|z|,u)^{1+\varepsilon}r_z,\,r_z\in (0,s).$

   By the covering theorem [20, с. 246] , from the covering $\{D(z,r_z)\}$ of the set $E$ it is possible to choose
   at most countable subcovering  $\{D(z_j,r_{z_j})\}$(in what follows we denote  $r_j:=r_{z_j}$ ),
   of multiplicity at most six with the following property:
\begin{equation}\mu (D(z_j,r_j))\ge B(|z_j|,u)^{1+\varepsilon}r_j.
\end{equation}
    We put $ R=r+(\log B(r,u))^{-1}$ and consider the sum
\begin{equation}\sum\limits_{r<|z_j|\le R}\mu (D(z_j,r_j))\ge
 B(r,u)^{1+\varepsilon} \sum\limits_{r<|z_j|\le R}r_j.
\end{equation}

    Above we applied (2.15) and monotonicity of the function  $B(r,u)$ .
    Applying Theorem
    Е and taking into account the radii of the exceptional disks,
    we obtain
    $$\sum\limits_{r<|z_j|\le R}\mu (D(z_j,r_j))\le
 6\mu \left(D\left(R+\frac{1}{B(r,u)^{1+\varepsilon}}\right)\right)
= 6\mu \left(D\left(r+\frac{1}{\log B(r,u)
}+\frac{1}{B(r,u)^{1+\varepsilon}}\right)\right)\le $$
\begin{equation}
\le 6\mu \left(D\left(r+\frac{2}{\log B(r,u) }\right)\right)\le6
C_{16} B(r,u)^{1+\delta},\,r\notin L \in FM.
\end{equation}

From (2.16 ) and (2.17) it follows (1.8). It remains to prove
 (2.14). From the condition $u(z)\le B(|z|,u)$ and the Jensen
 formula we deduce

     $$B(R,u)\ge
     \frac{1}{2\pi}\int\limits_0^{2\pi}u(Re^{i\varphi})\,d\varphi\ge
     \int\limits_r^R\frac{n(t,\mu_u)}{t}\,dt \ge
     n(r,\mu_u)\log\frac R r \ge C_{17}\frac{n(r,\mu_u)}{r\log B(r,u)}.
     $$
     Hence, by Theorem  E it follows that
\begin{equation}
n(r,\mu_u)\le C_{17}B(r,u)^{1+\delta/2}r\log B(r,u)<
B(r,u)^{1+\delta},\,r \notin L \in FM. \end{equation}
     From estimate (2.18) it follows (1.8), but only on the
     set
      $E_u:=\mathbb{C}\setminus N(\beta, s, \mu_u)$. In accordance with the remark after Theorem
      C, the entire function $f(z)$, the existence of which is proclaimed in that
      theorem, satisfies the inequality
\begin{equation}
\log |f(z)|\le u(z)+C_3(\varepsilon)(\log B(|z|,u)+\log |z|)\le 2
B(|z|,u),\,z\notin E_u.
 \end{equation}
     From estimate  (1.8) for the exceptional set $E_u$ it
     follows that for each
      $z$ there exists a sircumference  $S(0,t),\,t \in [r,R] $, where $r=|z|,R=r+(\log B(r,u))^{-1}$, all the points
       of which are
      $(\beta,s)$-normal with respect to the measure
      $\mu_u$. Applying the maximum principle for subharmonic functions and taking into account
    nondecrease of $B(r,u)$, from
      (2.19) we obtain the inequality
      $$
      \log|f(z)|\le 2B(R,u)+C_3(\varepsilon)\log |z| \le
      2B(r,u)^{1+\delta/2},
      $$
      which is satisfied for all $z$ but $|z|=r \notin L \in FM$. Again, by the Jensen formula
      we deduce
\begin{equation} n(r,\mu_{\log|f|})\le 2C_{17}B(R,u)^{1+\delta/2}r \log
B(r,u)< B(r,u)^{1+\delta},\,r\notin L \in FM.
 \end{equation}
       Combining (2.18) and (2.20), we come to (2.14).

        \textbf{The proof of Theorem  3.}
        \begin{lem} Let the conditions of theorem 3 be satisfied and $f(a)=0$. Then for each covering of the intersection
        $E\cap D(a,v(|a|)^{-1/2})$ by disks the sum of their radii is not less than  $v(|a|)^{-2C_4-1/2-\varepsilon}$.
        \end{lem}
  \textbf{The proof of Lemma 4} is similar to the proof of Lemma 3 and by this reason is omitted. We only  remark
  that here the properties of the function $v(R
  )$, proclaimed in Theorem 3, are applied.

   We now consider $b\in \mathbb{C}$ and
  the disk $D(b,r(b))$, where $r(b)=v(|b|)^{-1/2+\varepsilon}$.
  Either all the sircumferences $S(b,p r(b)),\,p\in [1/2,3/4],$ intersect $E$,
  or there exists a such sircumference, which does not intersect
    $E$. In the first case,
  the intersection $E\cap D(b,r(b)) $ is covered by disks having the sum of radii greater than
   $\frac 1 8 r(b)$. In the second case, there exists
  $a \in D(b,\frac 3 4 r(b))$ such that $f(a)=0$. Let us justify this statement.
  We suppose that $f$ has no zeros in the disk $D(b,p r(b))$.  We can write down the Poisson-Jensen
  formula for that disk
\begin{equation}
\frac{1}{2\pi}\int\limits_0^{2\pi}(\log|f(b+sr(b)e^{i\varphi})|-
v(|b+s r (b)e^{i\varphi}|))\,d\varphi=\int\limits_0^{p r (b)}
\frac{n(b,t,\log|f|)-n(b,t,v)}{t}\,dt.
\end{equation}

  The density of the Riesz measure  $\mu_v$ of the subharmonic function $v(|z|)$
  equals
   $\,\, \frac 1 {2\pi} \Delta v(|z|) \ge\,\,$\newline
    $\ge v(|z|)^{1-\varepsilon},\, z \in D(b,s r(b)) $.
   Therefore,
   \begin{equation}n(b,t,v)\ge v(|b|)^{1-\varepsilon}t^2.
\end{equation}
  By the assumption that the sircumference $S(b,s r(b))$ does not intersect the exceptional set
 it follows the left-hand side (2.21)
  is $O(\log v(|b|)),\,\,\,b \to \infty$, and from (2.22) it follows that the right-hand side
  (2.21) is greater than
  $ v(|b|)^{2\varepsilon}$ , i. e. we have a contradiction. Hence,
  the disk $D(b,s r(b))$ contains a zero $a$ of the function $f$. The disk
   $D(a,v(|a|)^{-1/2})\subset D(b,\frac 3 4 r(b))$ for
  sufficiently large   $|b|$ .
  We conclude that for every disk covering of the set $E \cap D(b, r(b))$ the sum of their radii is greater than
  $v(|a|)^{-2C_4-1/2-\varepsilon}\ge v(|b+r(b)|)^{-2C_4-1/2-\varepsilon}\ge ( v(|b|)^{-2C_4-1/2-
  \varepsilon})^{1+\varepsilon}\ge v(|b|)^{-2C_4-1/2-2\varepsilon}$.

  In the annulus $A(R,R+(\log v(R))^{-1})$ it is possible to put $\asymp R v(R)^{1-2\varepsilon}/\log v(R)$
  nonoverlapping disks of the form $D(b,r(b))$ (comparison of the areas).
 For every of them the sum of radii of the disks covering $E$ is
 greater than
   $v(R)^{-2C_4-1/2-2\varepsilon}$ (here and above we apply that by Theorem E in the annulus
    $A(R,R+(\log v(R))^{-1})$ the relation $r(b)\asymp v(R)^{-1/2+\varepsilon}$ holds) and multiplicity of the
    covering is bounded, therefore, the total sum of radii is greater than
   $v(R)^{-2C_4+1/2-4\varepsilon}$.

   \textbf{The proof of Theorem 4.} We hope that the qualified reader will accept a brief
   exposition.
       \begin{lem} Let the conditions of theorem 4 be satisfied and $f(a)=0$. Then for each disk covering of the intersection
        $E\cap D(a,(1-|a|)^{1+\rho/2})$ the sum of their radii is not less than
         \newline
         $(1-|a|)^{2C_5+1+\rho/2+\varepsilon}$.
        \end{lem}
        The proof of Lemma 5 is omitted, as it is similar to the
        proof of lemma 3.

              Next, by applying Lemma 5, we prove that for every
              disk covering of the intersection of the set
         $E$ with an arbitrary disk of the form  $D(b,(1-|b|)^{1+\rho/2-\varepsilon})$  the sum of their radii is not less than
       $(1-|b|)^{2C_5+1+\rho/2+\varepsilon}$.
         By the same arguments we deduce that the total sum of radii of the
         disk covering of the portion
         $E$ in the annulus $A(R,R+(1-R)/2)$ is not less than
        $\asymp\frac {1-R}{(1-R)^{2+\rho -2\varepsilon}}(1-R)^{2C_5+1+\rho/2+\varepsilon}=(1-R)^
        {2C_5-\rho/2+3\varepsilon}$.

        \textbf{The proof of Theorem 5.} We apply Theorem D, putting
         $\Omega = \mathbb{D},\,\beta=B(|z|,u)^{1+\varepsilon},\,s=\beta^{-1}$.
         by Theorem D under the conditions of Theorem 5
        there exists a function $f(z)$ analytic in $\mathbb{D}$
       satisfying (1.14) for $(\beta,\,s)$ -normal points with respect to the measures
         $\mu_u$
        and $\mu_{\log |f|}$. We denote $\nu=\mu_u+\mu_{\log|f|}$. In what follows the estimate
        \begin{equation}\nu(r)\le C_{18}B(r,u)^{1+\delta},\, r \notin L \in FM.
\end{equation} will be proved.
        Under the assumption that (2.23) is true, we estimate the
        size of the exceptional set
          $E$. By the definition of a normal point with respect to a measure,
         every point` $z\in E$ is the center of the disk $D(z,r_z)$ such that
\begin{equation}\nu(D(z,r_z))>
B(|z|,u)^{1+\varepsilon}r_z,\,r_z\in(0,B(|z|,u)^{-1-\varepsilon}).
\end{equation}

        By the covering theorem  [19, p.246], from the covering $\{D(z,r_z)\}$
        of the set $E$ we can choose at most countable subcovering
          $\{D(z_j,r_j)\}$(here $v_j:=r_{z_j}$) of multiplicity
          less than
          6. Applying estimate (2.24), we have
\begin{equation}
\sum\limits_{|z_j|\in[r,R]}\nu(D(z_j,r_j))\ge
B(r,u)^{1+\varepsilon} \sum\limits_{|z_j|\in[r,R]}r_j\,\,.
\end{equation}
        On the other hand, as it follows from  (2.23) and (2.24),
         $$
\sum\limits_{|z_j|\in[r,R]}\nu(D(z_j,r_j))\le 6
\nu\left(D\left(0,R+B(r,u)^{-1-\varepsilon}\right)\right)=$$\begin{equation}=
6 \nu\left(D\left(0,r+\frac {(1-r)^2}{\log B(r,u)
}+B(r,u)^{-1-\varepsilon}\right)\right) <
6\nu\left(D\left(0,r+\frac {2(1-r)^2}{\log B(r,u) }\right)\right)<
6 C_{18}B(r,u)^{1+\delta}.
\end{equation}
         Comparing (2.25) and (2.26), where we put $\delta=\varepsilon/2$,
         we obtain (1.15).

         It remains to prove (2.23). Without any restriction of generality,
         we can assume $u(0)=0$, then from the Poisson-Jensen
         formula it follows that
\begin{equation}B(R,u)\ge
     \frac{1}{2\pi}\int\limits_0^{2\pi}u(Re^{i\varphi})\,d\varphi\ge
     \int\limits_r^R\frac{n(t,\mu_u)}{t}\,dt \ge
     n(r,\mu_u)\log\frac R r \,.
\end{equation}
         From (2.27) and Theorem F we obtain that
$$ n(r,u)\log\left(1+\frac {(1-r)^2}{r\log B(r,u)}\right)\le
B(R,u) \le B(r,u)^{1+\delta},\, r \notin L \in FM , $$
        and, hence, the inequality
\begin{equation}n(r,u)\le B(r,u)^{1+\delta}\log B(r,u)/(1-r)^2\le
2B(r,u)^{1+\delta},\, r \notin L \in FM .
\end{equation}

         Hence, (1.15) holds on the set $E_u:=\mathbb{D}\setminus N(\beta,s,\mu_u)$. In accordance with the remark to
         Theorems C and D , the function $f$ satisfies the estimate from above
\begin{equation}\log |f(z)|\le u(z)+C_6(\varepsilon)\log
B(|z|,u)\le 2B(|z|,u),\,z\notin E_u.
\end{equation}

        From (2.28) and the proof of inequalities (2.25) and
        (2.26), it follows that for each point
          $z\in \mathbb{D}$ there exists a sircumference $S(0,T),\,T\in[(r+r)/2,R)$ ,
          all the points of which are $(\beta,s)-$normal with respect to the measure $\mu_u$.
          Then from
          (2.29) we obtain the inequality
\begin{equation}\log |f(z)|
\le 2B(T,u),\,z\in S(0,T).
\end{equation}
          Again, by the Poisson-Jensen formula we have
$$\int\limits_r^T \frac {\mu_{\log|f|}(t)} t\,dt\le\frac 1
{2\pi}\int\limits_0^{2\pi}\log|f(Te^{i\varphi})|\,d\varphi\le
2B(T,u), $$
           and from this we deduce
           $$n(r,\mu_{\log|f|})\log \frac T r \le 2B(T,u) \le
           2B(R,u)\le 2 B(r,u)^{1+\delta},\,r \notin L \in FM,
           $$
and, next,
\begin{equation}
           n(r,\mu_{\log|f|})\le \frac {2B(r,u)^{1+\delta}}{\log\left(1+\frac {(1-r)^2}{2\log
           B(r,u)}\right)}\le 4B(r,u)^{1+\delta}\log
           B(r,u)(1-r)^{-2}\le 4 B(r,u)^{1+2\delta}.
\end{equation}
           Therefore, by combining (2.31)
           with (2.28), inequality (2.23) is proved.

             \textbf{The proof of Theorem 6}, like to the proof of Theorem
              4, we expose briefly.
             We apply
               \begin{lem} Let the conditions of Theorem 6 be satisfied and $f(a)=0$. Then for every disk covering
        of the intersection $E\cap D(a,v(|a|)^{-1/2})$ the sum of their radii is not less than
         $v(|a|)^{-2C_7-1/2-\varepsilon}.$
        \end{lem}
       Next, we show that for every disk covering of the intersection of the set
        $E$ with an arbitrary disk  $D(b,v(|b|)^{-1/2+\varepsilon})$ the sum of the radii of the covering disks
      is not less than $v(|b|)^{-2C_7-1/2-\varepsilon}$.
         In the annulus $A(R,R+(1-R)^2/\log v(R)))$ we can arrange about  $\frac {(1-R)^2 v(R)^{1-2\varepsilon}}{\log v(R)}$
  of nonoverlapping disks of such a form (comparison of the areas). In every of them the sum of
  radii of the disks covering of the set
    $E$, is greater than
   $v(R)^{-2C_7-1/2-\varepsilon}$ (In this place and above we use that by Theorem
   F in the annulus $A(R,R+(1-R)^2/\log v(R))$ the relation $r(b)\asymp
   v(R)^{-1/2+\varepsilon}$holds) and multiplicity of the covering
   is finite, hence, the total sum of radii is greater than
    $v(R)^{-2C_7+1/2-31\varepsilon}$.

\end{document}